\documentclass[12pt]{article}
\usepackage{graphicx}
\usepackage{amssymb}
\usepackage{amsmath}
\numberwithin{equation}{section}
\def\be{\begin{equation}\label}
\def\e{\end{equation}}

\author{Lev Sakhnovich}

\title{About the Rationality of the Knizhnik-Zamolodchikov Equation Solution }

\date{}

\begin{document}

\maketitle

735 Crawford ave., Brooklyn, 11223, New York, USA.\\
 E-mail address: lev.sakhnovich@verizon.net

\begin{abstract}
 In the paper the solution of KZ system (n=4, m=2) is constructed in the explicit
 form in  terms of the hypergeometric functions. We proved that the corresponding
 solution is rational when the parameter $\rho$ is integer. We show that in the case
 (n=4, m=5, $\rho$ is integer) the corresponding KZ system hasn't got a rational solution.
\end{abstract}

\textbf{Mathematics Subject Classification (2000).} Primary 35C05; Secondary 20B05, 15A18.\\
 \textbf{Keywords.} Symmetric group, Natural representation, Young tableau, Integer eigenvalues.
 \section{Introduction}
 1. We consider the differential system
\begin{equation}
\frac{{\partial}W(z_{1},z_{2},...,z_{n})}{{\partial}z_{k}}={\rho}A_{k}(z_{1},z_{2},...z_{n})W,
\quad  1{\leq}j{\leq}n,\end{equation} where
 $A_{k}(z_{1},z_{2},...,z_{n})$ and $W(z_{1},z_{2},...,z_{n})$ are $m{\times}m$ matrix
functions. We suppose that $A_{k}(z_{1},z_{2},...,z_{n})$ has the
form
\begin{equation}
A_{k}(z_{1},z_{2},...,z_{n})=\sum_{j=1,
j{\ne}k}^{n}\frac{P_{k,j}}{z_{k}-z_{j}},\end{equation} where
$z_{k}{\ne}z_{j}$ if $k{\ne}j$. Here the matrices $P_{j,k}$
are connected with the matrix representation of the symmetric group
$S_{n}$ and are defined by formulas (2.1)-(2.4). We note that the
well-known Knizhnik- Zamolodchikov equation has the form (1.1), (1.2)
(see [4]). This system has found applications in several  areas of
mathematics and physics (see [4],[12]). In
paper ([11],section 2) we prove the following assertion:\\
\textbf{Theorem 1.1.} \emph{The fundamental solution of KZ system (1.1),(1.2) is
rational,when ${\rho}$ is integer and matrices $P_{k,j}$ give
the natural representation of symmetric group $S_{n}$.}\\
We note that the natural representation is the sum of the 1-representation and the
 irreducible (n-1)-representation. We name the corresponding irreducible (n-1)-representation \emph{the  natural representation} as well.\\
 In the present  paper the solution of KZ system (n=4, m=2) is constructed in the explicit form in the terms of the hypergeometric functions. We proved that the corresponding solution is rational,if $\rho$ is integer.
 A number of authors (see [3],[5])  affirmed that the solution of
KZ equation is rational for all the representation of symmetric group $S_{n}$, if
the parameter $\rho$ is integer. In this paper we construct the example  (n=4 and n=5) for which ${\rho}$ is integer but the corresponding solution is irrational. This result leads to the following open problem.\\
\textbf{Open Problem 1.1} \emph{Let $\rho$ be integer and let $P_{k,j}$ be an irreducible representation of symmetric group $S_{n}$. To find the conditions under which  the corresponding fundamental solution is rational.}\\Here can be useful the following necessary condition
{see [9]):\\
\textbf{Proposition 1.1} \emph{If $\rho$ is integer and the fundamental solution of  system (1.1), (1.2) is rational then  all the eigenvalues of matrices}
\begin{equation}Q_{k}=\sum_{j{\ne}k,j=1}^{n}P_{k,j},\quad 1{\leq}k{\leq}n,
\end{equation}\emph{ are integer.}\\
In the before mentioned  case (n=4, m=2) the formulated necessary condition is  fulfilled.
\section{Knizhnik-Zamolodchikov Equation Solution, \\ Case n=4, m=2.}
 In the case n=4, m=2 we have the following matrix irreducible representation (see [2], [6]):
\begin{equation} P_{1,2}=P_{3,4}=\left(
                           \begin{array}{cc}
                             1 & 0 \\
                             0 & -1 \\
                           \end{array}
                         \right), \quad P_{2,4}= P_{1,3}=\left(
                           \begin{array}{cc}
                             -1/2 & -\sqrt{3}/2 \\
                            -\sqrt{3}/2  & 1/2 \\
                           \end{array}
                         \right), \end{equation} \begin{equation}P_{2,3}= P_{1,4}=\left(
                           \begin{array}{cc}
                             -1/2 & \sqrt{3}/2 \\
                            \sqrt{3}/2  & 1/2 \\
                           \end{array}
                         \right) .\end{equation}
Following A.Varchenko [12] we change the variables
\begin{equation}u_{1}=z_{1}-z_{2},\quad
u_{k}=\frac{z_{k}-z_{k+1}}{z_{k-1}-z_{k}},\quad
(2{\leq}k{\leq}3),\end{equation} \begin{equation}
u_{4}=z_{1}+z_{2}+z_{3}+z_{4}.\end{equation} The KZ system takes the
following form
\begin{equation}\frac{{\partial}W}{{\partial}u_{j}}=
{\rho}H_{j}(u)W,\quad (1{\leq}j{\leq}4),\end{equation}where
$u=(u_{1},u_{2},u_{3},u_{4}).$
  We have [11]:
\begin{equation}\frac{{\partial}W}{{\partial}u_1}=\rho{\left(\frac{\Omega_1}{u_1}\right)}W,\end{equation}
\begin{equation}\frac{{\partial}W}{{\partial}u_2}=\rho{\left(\frac{\Omega_{2}}{u_2}+
\frac{P_{1,3}}{1+u_2}+\frac{P_{1,4}(1+u_{3})}{1+u_2+u_{2}u_{3}}\right)}W,\end{equation}
\begin{equation}\frac{{\partial}W}{{\partial}u_3}=\rho{\left(\frac{P_{4,3}}{u_3}+
\frac{P_{4,2}}{1+u_3}+\frac{P_{4,1}u_{2}}{1+u_2+u_{2}u_{3}}\right)}W,\end{equation}where                         
$H_{4}(u)=0$ and
\begin{equation}P_{r}=\sum_{j>r}P_{j,r},\quad \Omega_{s}=P_{s}+P_{s+1}+...+P_{4}.\end{equation}
It follows from (2.1),(2.2) and (2.9)  that
$\Omega_{1}=\Omega_{2}=0.$ Hence KZ system (2.7), (2.8) takes the form
\begin{equation}\frac{{\partial}W}{{\partial}u_2}=\rho{\left(
\frac{P_{1,3}}{1+u_2}+\frac{P_{1,4}(1+u_{3})}{1+u_2+u_{2}u_{3}}\right)}W,\end{equation}
\begin{equation}\frac{{\partial}W}{{\partial}u_3}=\rho{\left(\frac{P_{4,3}}{u_3}+
\frac{P_{4,2}}{1+u_3}+\frac{P_{4,1}u_{2}}{1+u_2+u_{2}u_{3}}\right)}W.\end{equation}
We introduce the
matrix function
\begin{equation}F(y)=W(y)(1+y)^{\rho}(1+y+yz)^{\rho},\quad
y=u_{2},\quad u_{3}=z.\end{equation}Relations (2.10) and (2.12) imply that
\begin{equation}\frac{dF}{dy}=\rho{\left(\frac{P_{1,3}+I}{1+y}+
\frac{(P_{1,4}+I)(1+z)}{1+y+yz}\right)}F.\end{equation}We consider the
constant vectors
\begin{equation}w_{1}=\mathrm{col}[\sqrt{3},1],\quad
w_{2}=\mathrm{col}[1,\sqrt{3}].\end{equation}It is easy to see that
\begin{equation}(P_{1,3}+I)w_{1}=0,\quad
(P_{1,3}+I)w_{2}=-\sqrt{3}w_{1}+2w_{2},\end{equation}
\begin{equation}(P_{1,4}+I)w_{1}=\sqrt{3}w_{2},\quad
(P_{1,4}+I)w_{2}=2w_{2},\end{equation} We represent $F(y)$ in the
form
\begin{equation}F(y)=\phi_{1}(y)w_{1}+\phi_{2}(y)w_{2}\end{equation}
and substitute it in (2.13). In view of (2.15) and (2.16) we have
\begin{equation}\phi_{1}^{\prime}(y)=-\frac{\sqrt{3}\rho}{1+y}\phi_{2}(y),\end{equation}
\begin{equation}
\phi_{2}^{\prime}(y)=\rho\left(\frac{2\phi_{2}(y)}{1+y}+\frac{(2\phi_{2}(y)+\sqrt{3}\phi_{1})(1+z)}{1+y+yz}\right).
\end{equation}
\begin{equation}(1+y)(\frac{1}{z+1}+y)\phi_{1}^{\prime\prime}(y)+B(z,y)\phi_{1}^{\prime}(y)
+3{\rho}^{2}\phi_{1}(y)=0,\end{equation}where
\begin{equation}B(z,y)=\frac{1}{z+1}+y-2\rho[2(\frac{1}{z+1}+y)+\frac{z}{z+1}]
\end{equation}
Changing the variable $y=\frac{-z}{z+1}v-1,\quad \phi(v)=\phi_{1}(\frac{-z}{z+1}v-1)$
we reduce equation (2.20) to the following form
\begin{equation}v(1+v)\phi^{\prime\prime}(v)+[1+v-2\rho(1+2v)]\phi^{\prime}(v)
+3{\rho}^{2}\phi(v)=0.\end{equation} By introducing $\psi(v)=\phi(-v)$ we obtain   Gauss
hypergeometric equation [1]:
\begin{equation}v(1-v)\psi^{\prime\prime}(v)+[\gamma-(\alpha+\beta+1)y]\psi^{\prime}(v)
-\alpha\beta\psi(v)=0,\end{equation}where
\begin{equation}\alpha=-\rho,\quad \beta=-3\rho,\quad
\gamma=1-2\rho.\end{equation}In our paper [11] we have proved the following assertion.\\
\textbf{Proposition 2.1.} \emph{The solutions of Gauss
hypergeometric equation (2.23) are rational functions if
$$\alpha=-\rho,\quad \beta=-3\rho,\quad
\gamma=1-2\rho,$$ where $\rho$ is integer.}\\
Let $\psi_{1}(v)$ and $\psi_{2}(v)$ be
linearly independent rational solutions of equation (2.22). Then the vector
functions
\begin{equation}Y_{k}(y,z)=(1+y)^{-\rho}(1+y+yz)^{-\rho}X_{k}(y,z),\quad (k=1,2),\end{equation} where
\begin{equation}X_{k}(y,z)=\psi_{k}\left(\frac{(y+1)(z+1)}{z}\right)w_{1}+
\frac{(y+1)(z+1)}{z\sqrt{3}\rho}\psi_{k}^{\prime}\left(\frac{(y+1)(z+1)}{z}\right)w_{2}
\end{equation}
 are the
solutions of system (2.13).
 Hence we deduced the assertion\\
\textbf{Proposition 2.2.} \emph{The fundamental solution $W_{1}(y,z,\rho)$ of
system (2.10) is defined by relation}
\begin{equation}W_{1}(y,z,\rho)=[Y_{1}(y,z,\rho),Y_{2}(y,z,\rho)],\end{equation}
\emph{where} $y=u_{2},\quad z=u_{3}.$\\
\textbf{Proposition 2.3.} \emph{If $\rho=-1$ then the solutions $Y_{1}(y,z,-1)$ and $Y_{2}(y,z,-1)$ of equation (2.10) have the forms}
\begin{equation}Y_{1}(y,z,-1)=-z(y+1)w_{1}+\frac{z(z+1)(y+1)^{2}}{\sqrt{3}(1+y+yz)}w_{2},\end{equation}
\begin{eqnarray}&& \nonumber Y_{2}(y,z,-1)=[\frac{z^{2}(1+y+yz)}{(z+1)^{2}(y+1)}+
\frac{z(1+y+yz)}{(z+1)}]w_{1}
\\ &&
-[\frac{2z^{2}(1+y+yz)}{\sqrt{3}(z+1)^{2}(y+1)}+
\frac{z(1+y+yz)}{\sqrt{3}(z+1)}]w_{2}.\end{eqnarray}
\section{Common solution, Consistency}
We note that the KZ system (2.10), (2.11) is consistent [11].
In section 2 we have considered only the first  equation of  system (2.10), (2.11).
Now using this result we shall construct the common solution of the  KZ system (2.10),(2.11).
To do it  let us consider  equation (2.11) in case when $u_{2}=0,\quad u_{3}=z.$ We have
\begin{equation}\frac{{\partial}W_{2}}{{\partial}z}=\rho{\left(\frac{P_{4,3}}{z}+
\frac{P_{4,2}}{1+z}\right)}W_{2}.\end{equation}
This equation can be solved in the same way as (2.10).
 We introduce the
matrix function
\begin{equation}G(z)=W_{2}(z)z^{\rho}(1+z)^{\rho},\quad
z=u_{3}.\end{equation}Relations (3.1) and (3.2) imply that
\begin{equation}\frac{dG}{dz}=\rho{\left(\frac{P_{4,3}+I}{z}+
\frac{P_{4,2}+I}{1+z}\right)}G.\end{equation}We consider the
constant vectors
\begin{equation}v_{1}=\mathrm{col}[0,2],\quad
v_{2}=\mathrm{col}[1,-\sqrt{3}].\end{equation}It is easy to see that
\begin{equation}(P_{4,3}+I)v_{1}=0,\quad
(P_{4,2}+I)w_{1}=-\sqrt{3}v_{2},\end{equation}
\begin{equation}(P_{4,3}+I)v_{2}=\sqrt{3}v_{1}+2v_{2},\quad
(P_{4,2}+I)v_{2}=2v_{2},\end{equation} We represent $G(z)$ in the
form
\begin{equation}G(z)=\phi_{1}(z)v_{1}+\phi_{2}(z)v_{2}\end{equation}
and substitute it in (3.3). In view of (3.5) and (3.6) we have
\begin{equation}\phi_{1}^{\prime}(z)=\frac{\sqrt{3}\rho}{z}\phi_{2}(z),\quad
\phi_{2}^{\prime}(z)=\rho(\frac{2\phi_{2}(z)}{z}+\frac{2\phi_{2}(z)}{1+z}-\frac{\sqrt{3}\phi_{1}(z)}{1+z}).
\end{equation}It follows from (3.8) that
\begin{equation}z(1+z)\phi_{1}^{\prime\prime}(z)+[1+z-2\rho(1+2z)]\phi_{1}^{\prime}(z)
+3{\rho}^{2}\phi_{1}(z)=0.\end{equation}By introducing
$\psi(z)=\phi_{1}(-z)$ we reduce equation (3.9) to Gauss
hypergeometric equation (2.23), (2.24).
Let $\psi_{1}(z)$ and $\psi_{2}(z)$ be
linearly independent solutions of equation (2.23). Then the vector
functions
\begin{equation}U_{k}(z,\rho)=[\psi_{k}(-z)v_{1}-\frac{z}{\sqrt{3} \rho}\psi_{k}^{\prime}(-z)v_{2}]
z^{-\rho}(1+z)^{-\rho},\quad k=1,2\end{equation} are the
solutions of system (3.1). Hence we deduced the assertion\\
\textbf{Proposition 3.1.} \emph{The fundamental solution $W_{2}(z,\rho)$ of
system (3.1) is defined by the relation}
\begin{equation}W_{2}(z,\rho)=[U_{1}(z,\rho),U_{2}(z,\rho)].\end{equation}
Let us consider separately the case when $\rho=-1.$ In this case equation (3.9) takes the form
\begin{equation}z(1+z)\phi_{1}^{\prime\prime}(z)+(3+5z)\phi_{1}^{\prime}(z)
+3\phi_{1}(z)=0.\end{equation}
It is easy to check by direct calculation that
the functions
\begin{equation}\phi_{1,1}=\frac{1}{1+z},\quad \phi_{1,2}=\frac{1-z}{z^{2}}\end{equation}
are solutions of equation (3.12). Hence the following assertion is true:\\
\textbf{Proposition 3.2.} \emph{If $\rho=-1$ then the solutions $U_{1}(z,-1)$ and $U_{2}(z,-1)$ of equation (3.12) have the forms}
\begin{equation}U_{1}(z,-1)=zv_{1}+\frac{z^{2}}{\sqrt{3}(1+z)}v_{2},\quad
U_{2}(z,-1)=\frac{1-z^{2}}{z}v_{1}-\frac{(z-2)(z+1)}{\sqrt{3}z}v_{2}\end{equation}
We have constructed the fundamental solutions $W_{1}(y,z,\rho)$ and $W_{2}(y,z,\rho)$
of equations (2.10) and (3.1) respectively. Now using our previous results (see[],section 3)
we obtain the main theorem of the present paper.\\
\textbf{Theorem 3.1.}\emph{The $2{\times}2$ matrix function }
\begin{equation}W(y,z,\rho)=W_{1}(y,z,\rho)W_{1}^{-1}(0,z,\rho)W_{2}(0,z,\rho)\end{equation}
\emph{is the fundamental  solution
of the KZ system (2.10),(2.11), where $y=u_{2},\quad z=u_{3}$. If $\rho$ is integer this fundamental solution is rational.}
\section{System with Non-rational Solution.}
Let us consider the case n=5, m=5.\\ In this case we have the following matrix irreducible representation (see [2], [6]):
\begin{equation} P_{1,2}=\left(
 \begin{array}{ccccc}
                             1 & 0 & 0 & 0 & 0 \\
                             0 & 1 & 0 & 0 & 0 \\
                             0 & 0 & -1 & 0 & 0  \\
                             0 & 0 & 0 &  1 & 0 \\
                             0 & 0 & 0 & 0 & -1
                           \end{array}\right)
                       , \end{equation}
           \begin{equation}                P_{1,3}=\left(
                           \begin{array}{ccccc}
                             1 & 0 & 0 & 0 & 0 \\
                             0 & -1/2 &-\sqrt{3}/2 & 0 & 0 \\
                             0 & -\sqrt{3}/2 & 1/2 & 0 & 0 \\
                             0 & 0 & 0 & -1/2 & -\sqrt{3}/2 \\
                             0 & 0 & 0 & -\sqrt{3}/2 & 1/2
                           \end{array}
                         \right) ,\end{equation}

\begin{equation}P(1,4)=\left(\begin{array}{ccccc}
                               -1/3 & -\sqrt{2}/3 & -\sqrt{6}/3 & 0 & 0 \\
                               -\sqrt{2}/3 & 5/6 & -\sqrt{3}/6 & 0 & 0 \\
                               -\sqrt{6}/3  & -\sqrt{3}/6  & 1/2 & 0 & 0 \\
                               0 & 0 & 0 & -1/2 & \sqrt{3}/2\\
                               0 & 0 & 0 & \sqrt{3}/2 & 1/2
                             \end{array}\right),\end{equation}
\begin{equation} P(1,5)=\left(\begin{array}{ccccc}
                                -1/3 & \sqrt{2}/9 & \sqrt{6}/9 & -4/9 & -4\sqrt{3}/9 \\
                                \sqrt{2}/9 & -19/54 & 23\sqrt{3}/54 & -8\sqrt{2}/27 & 4\sqrt{6}/27 \\
                                \sqrt{6}/9 &  23\sqrt{3}/54 & 1/2 & 4\sqrt{6}/27 & 0 \\
                                -4/9 & -8\sqrt{2}/27 & 4\sqrt{6}/27 & 37/54 & -5\sqrt{3}/54 \\
                                -4\sqrt{3}/9 & 4\sqrt{6}/27 & 0 & -5\sqrt{3}/54 & 1/2
                              \end{array}\right).\end{equation}Using the previous representations of $P(1,k)$
                              we deduce that
the matrix $Q_{1}=P(1,2)+P(1,3)+P(1,4)+P(1,5)$ has the form
\begin{eqnarray}&\nonumber
T_{-1}=\Big(C_1 \quad C_2  \quad \ldots \quad C_5 \Big), \quad &
C_1=\left(\begin{array}{c}4/3  \\ -2\sqrt{2}/9 \\-(2\sqrt{2/3})/9\\-4/9 \\-4/(3\sqrt{3}) 
\end{array}\right), \\ \nonumber
&
C_2=\left(\begin{array}{c}  -2\sqrt{2}/9\\53/54\\7/(9\sqrt{3})-\sqrt{3}/2 \\-8\sqrt{2}/27\\(4\sqrt{2/3})/9
\end{array}\right), \quad &
C_3=\left(\begin{array}{c}  -(2\sqrt{2/3})/3\\7/(9\sqrt{3})-\sqrt{3}/2\\1/2 \\ (4\sqrt{2/3})/9 \\0
\end{array}\right), 
\\& \nonumber
C_4=\left(\begin{array}{c}  -4/9\\-8\sqrt{2}/27\\(4\sqrt{2/3})/9\\37/54\\-5/(18\sqrt{3})
\end{array}\right), \quad  &
C_5=\left(\begin{array}{c}  -4/(3\sqrt{3})\\(4\sqrt{2/3})/9\\0\\-5/(18\sqrt{3})\\ 1/2
\end{array}\right).
\end{eqnarray}

The eigenvalues of the matrix $T_{-1}$ are defined by the relations
\begin{equation} \lambda_{1,2}=(17{\pm}\sqrt{433})/18,\quad \lambda_{3}=5/3, \quad \lambda_{4}=1/3,\quad \lambda_{5}=1/9.\end{equation}
We see that  all the eigenvalues of $Q_{1}$ are non-integer.
Hence according to Proposition 1.1 we obtain the assertion.\\
\textbf{Proposition 4.1} \emph{Let $\rho$ be integer.The KZ system in the case (n=5,m=5) has no rational solutions.}

\end{document}